\documentclass[3p,times]{elsarticle}

\usepackage{amssymb}
\usepackage{amsmath}

\DeclareMathOperator{\e}{e}

\newtheorem{lemma}{Lemma}

\newtheorem{Example}{Example}
\newenvironment{example}{\begin{Example}\normalfont}{\end{Example}}

\begin{document}

\begin{frontmatter}




\title{An efficient numerical approach for delayed logistic models}

\author{Josef Rebenda\corref{auth1}}
\author{Zden\v{e}k \v{S}marda\fnref{auth2}}
\fntext[auth1]{CEITEC BUT, Brno University of Technology, Purkynova 123, 61200 Brno, Czech Republic}
\cortext[auth2]{CEITEC BUT, Brno University of Technology, Purkynova 123, 61200 Brno, Czech Republic}

\address{}

\begin{abstract}
 
In the paper an efficient semi-analytical approach based on the method of steps and differential transformation is proposed for numerical approximation of solutions of retarded logistic models of delayed and neutral type, including models with several constant delays. Algorithms for both commensurate and non-commensurate delays are described, applications are shown in examples. Validity and efficiency of the presented algorithms is compared with variational iteration method, Adomian decomposition method and polynomial least squares method numerically. Matlab package DDE23 is used to produce reference numerical values.

\end{abstract}


\end{frontmatter}


\section{Introduction}
The logistic equation (sometimes called the Verhulst model or logistic growth curve) is a
model of population growth first published by Pierre-Francois Verhulst (1845). The model is continuous in time, 
but a modification of the continuous equation to a discrete
quadratic recurrence equation is also widely studied. The continuous version of the logistic model has the form
\begin{equation}\label{l1}
N'(t) = rN(t)\left( 1-\frac{N(t)}{K}\right),
\end{equation}
where $r>0$ is the intrinsic growth rate and $K>0$ is the carrying capacity of the population and $N(t)$ is the size of the population at time $t$.
To make the logistic equation more realistic, Hutchinson \cite{hutch} proposed incorporating the
effect of delay and he introduced the delayed logistic equation
\begin{equation}\label{l2}
N'(t) = rN(t)\left( 1-\frac{N(t-\tau)}{K}\right),
\end{equation}
where $\tau >0$ is a time delay which represents maturation time. Hutchinson suggested that the equation \eqref{l2} can be used to model the dynamics of a single species
population growing towards a saturation level $K$ with a constant reproduction rate $r$ (see \cite{gopal}, \cite{kuang}).
Another form, studied by Gy\"ori \cite{gyori}, is the equation
with $n$ delays defined for $t \geq 0$ by
\begin{equation}\label{l3}
N'(t) = N(t)\left(a- \sum_{j=1}^n b_j N(t-\tau_j)\right).
\end{equation}
This equation is supposed to describe a situation where several of the processes affecting the population occur with
different time delays.

Delay logistic type equations  have also found applications in ecology (\cite{gopal}, \cite{pielou}),  where the following neutral delay model \cite{pepe} 
for forest  evolution based on Pearl-Vehulst equation has been proposed:
\begin{equation}\label{l4}
N'(t) = rN(t)\left( 1-\frac{N(t-\tau)+ c N'(t-\tau)}{K}\right),
\end{equation}
where $x,r,K$ are the tree population, the intrinsic growth rate and the enviromental carrying capacity. The additional term $c N'(t-\tau)$
to the usual logistic equation  is introduced for soil depletion and erosion.

Oscillatory properties of solutions of another form of neutral logistic equation were studied by Gy\"ori and Ladas \cite{gyori2}:
\begin{equation}\label{l5}
N'(t) = N(t)\left\{ r \left[1-\frac{N(t-\sigma)}{K} \right]+ c \frac{N'(t-\tau)}{N(t-\tau)}\right\},
\end{equation}
where the term $ c \frac{N'(t-\tau)}{N(t-\tau)}$ is associated with per capita growth rate of the population at time $t-\tau$.

There are a lot of papers devoted to investigation of qualitative properties of delayed logistic models. On the other hand, semianalytical methods convenient for solving of delayed models are in forefront of the study in last two decades. For all we mention papers \cite{dehghan}, \cite{caruntu} and \cite{khan}, where these models are treated numerically. However, the calculations and results are complicated, therefore we propose an easy applicable approach in this paper.

\section{Preliminaries}
\subsection{Functional Differential Equation}
We will consider the following general functional differential equation of $n$-th order and then apply the main results to a special case of delayed or neutral logistic equation:
\begin{equation}\label{1}
u^{(n)}(t) = f(t,u(t),u'(t),\dots,u^{(n-1)}(t),\mathbf{u}_1(t-\tau_1 ),\mathbf{u}_2(t-\tau_2 ), \dots, \mathbf{u}_r(t-\tau_r ) ),
\end{equation}
where $\mathbf{u}_i(t-\tau_i )= (u(t-\tau_i),u'(t-\tau_i),\dots,u^{(m_i)}(t-\tau_i))$ is $m_i$-dimensional vector function, $m_i \leq n$, $i=1,2,\dots,r$, $r \in \mathbb N$ and $f \colon [t_0,\infty) \times R^n \times R^\omega$ is a continuous function, where $\omega = \sum\limits_{i=1}^r m_i$.\\
Let $t^*= \max\{\tau_1,\tau_2,\dots,\tau_r\}$, $m= \max\{m_1,m_2,\dots,m_r\}$, $m \leq n $.  In case $m=n$ equation \eqref{1} is of neutral type, otherwise it is delayed differential equation. Initial function $\phi(t)$ needs to be assigned  for equation \eqref{1} on the interval $[t_0-t^*, t_0]$.
Furthermore, for the sake of simplicity, we assume that $\phi(t) \in C^n([t_0-t^*,t_0])$.

Various methods such as homotopy analysis method (HAM) \cite{khan}, \cite{wang}, \cite{alomari}, homotopy  perturbation method (HPM) \cite{shakeri}, variational iteration method (VIM) \cite{chen},  Adomian decomposition method (ADM)  \cite{evans}, \cite{cocom}, Taylor polynomial method \cite{sezer}, Taylor collocation method \cite{bellour}
and differential transformation method (DTM) \cite{arikoglu}, \cite{karakoc}, have been considered to approximate solutions of certain classes of equation \eqref{1} in a series form. 

The crucial idea of our concept is to combine differential transformation method and general method of steps (more details on method of steps can be found for instance in monographies Kolmanovskii and Myshkis \cite{kolmanovski}, Hale and Verduyn Lunel \cite{hale} or Bellen and Zennaro \cite{bellen}).
This approach enables us to replace the terms involving delay with initial function and its derivatives. Consequently, the original Cauchy problem for delayed or neutral differential equation is reduced to Cauchy problem for ordinary differential equation. Further, while ADM, HAM, HPM and VIM require initial guess approximation  and symbolic computation of necessary derivatives and, in general, $n$-dimensional integrals in iterative schemes, presented method is different: Cauchy problem for FDE is reduced to a system of recurrence algebraic relations.

\subsection{Overview of Differential Transformation}

Differential transformation of a real function $u(t)$ at a point $t_0 \in \mathbb R$ is $\mathcal D \{ u(t) \} [t_0] = \{ U(k) \}_{k=0}^{\infty}$,
where $U(k)$, differential transformation of the $k-$th derivative of function $u(t)$ at $t_0$, is defined as
\begin{equation}\label{2}
U(k) = \frac{1}{k!} \left[ \frac{d^ku(t)}{dt^k} \right]_{t=t_0},
\end{equation}
provided that the original function $u(t)$ is sufficiently smooth. Inverse differential transformation of $\{ U(k) \}_{k=0}^{\infty}$ is defined as follows:
\begin{equation}\label{3}
u(t) = \mathcal D^{-1} \Bigl\{ \{ U(k) \}_{k=0}^{\infty} \Bigr\} [t_0]= \sum_{k=0}^{\infty}U(k)(t-t_0)^k. 
\end{equation}
As we can observe from \eqref{3}, this semi-analytical numerical technique is based on Taylor series.
In real applications the function $u(t)$ is expressed by a finite series 
\begin{equation}\label{3f}
u(t) = \sum_{k=0}^{N}U(k)(t-t_0)^k. 
\end{equation}
Plenty of transformation formulas can be derived from definitions \eqref{2} and \eqref{3}, we recall the following which we will use later in numerical experiments:
\begin{lemma}\label{lemma1}
Assume that $F(k)$, $G(k)$, $H(k)$ and $U_i(k)$, $i=1,\dots,n$, are differential transformations of functions $f(t)$, $g(t)$, $h(t)$ and
$u_i(t)$, $i=1,\dots,n$, respectively. Then
\begin{gather*}
\begin{array}{lllcl}
i) & {\rm  If} & f(t) = {\displaystyle \frac{d^ng(t)}{dt^n}}, & {\rm then} & F(k) = {\displaystyle \frac{(k+n)!}{k!}}G(k+n). \\[2mm]
ii)& {\rm  If} & f(t) = g(t)h(t), & {\rm then} & F(k) = \sum_{l=0}^k G(l)H(k-l). \\[2mm]
iii)& {\rm  If} & f(t) = t^n, & {\rm then} & F(k) =\delta (k-n),\ \delta\ {\rm is\ the\ Kronecker\ delta\  symbol}, t_0=0. \\[1mm]
iv)& {\rm  If} & f(t) = e^{\lambda t}, & {\rm then} & F(k) = {\displaystyle \frac{\e^{a \lambda} \lambda^k}{k!}}, t_0=a. \\[2mm]
v)& {\rm  If} & f(t) = t^n, & {\rm then} & F(k) =\begin{cases}
  \binom{n}{k} a^{n-k} & \text{if } k \leq n,\\
  0 & \text{if } k>n 
 \end{cases}, t_0=a, n \in \mathbb{N}_0. \\[2mm]
vi) &{\rm  If} & f(t) = \prod_{i=1}^n u_i(t),& {\rm then}&
\end{array}
\\
F(k)= \sum_{s_1=0}^k \sum_{s_2=0}^{k-s_1} \dots \sum_{s_n=0}^{k-s_1-\cdots s_{n-1}} U_1(s_1)\dots U_{n-1}(s_{n-1}) U_n(k-s_1-\dots -s_n).
\end{gather*}
\end{lemma}

{\bf Remark.} Transformation formulas for shifted arguments also can be found in the literature. However, this approach is not convenient for solving initial problems for retarded differential equations. A convenient approach is presented in the following section.





\section{Main Results}\label{sMain}
Recall equation \eqref{1} 
\begin{equation}
u^{(n)}(t) = f(t,u(t),u'(t),\dots,u^{(n-1)}(t),\mathbf{u}_1(t-\tau_1 ),\mathbf{u}_2(t-\tau_2 ), \dots, \mathbf{u}_r(t-\tau_r ) ),
\tag{\ref{1}}
\end{equation}
and consider it  subject to initials conditions
\begin{equation}\label{6}
u(t_0)=u_{1,0}, u'(t_0)=u_{1,1}, \dots, u^{(n-1)}(t_0) = u_{1,n-1}
\end{equation}
and subject to initial function $\phi(t)$ on interval $[t_0-t^*, t_0]$ such that
\begin{equation}\label{6'}
\phi(t_0) = u(t_0), \phi'(t_0)= u'(t_0), \dots, \phi^{(n-1)}(t_0)= u^{(n-1)}(t_0).
\end{equation}
First we apply the method of steps. We substitute the initial function $\phi(t)$ and its derivatives in all places where unknown functions with deviating arguments and derivatives of that functions appear. Then equation \eqref{1} changes to ordinary differential equation
\begin{equation}\label{7}
u_1^{(n)}(t) = f(t,u_1(t),u_1'(t),\dots,u_1^{(n-1)}(t),\mathbf{\Phi}_1(t-\tau_1),\mathbf{\Phi}_2(t-\tau_2), \dots, \mathbf{\Phi}_r
(t-\tau_r) ),
\end{equation}
where $\mathbf{\Phi}_i(t-\tau_i )= (\phi(t-\tau_i),\phi'(t-\tau_i),\dots,\phi^{(m_i)}(t-\tau_i))$, $m_i \leq n$, $i=1,2,\dots,r$.
Now applying DTM we get recurrence  equation
\begin{equation}\label{8}
\frac{(k+n)!}{k!} U_1(k+n) = \mathcal F_1 \Bigl( k, U_1(k), U_1(k+1), \dots, U_1(k+n-1)\Bigr),
\end{equation}
where
$$
U_1(0)= u(t_0)= u_{1,0},\ U_1(1)= u'(t_0) =u_{1,1}, \dots , U_1(n-1)= \frac{1}{(n-1)!} u^{(n-1)}(t_0) = \frac{u_{1,n-1}}{(n-1)!}.
$$

Using transformed initial conditions and then inverse transformation rule, we obtain a solution of
equation \eqref{1} in the form of Taylor series
$$
u_1(t) = \sum_{k=0}^{\infty} U_1(k)(t-t_0)^k
$$
on the interval $[t_0, t_0+ \alpha]$, where $\alpha = \min\{\tau_1,\tau_2, \dots, \tau_r\}$, and $ u_1(t) = \phi(t)$ on the interval $[t_0-t^*,t_0]$.

The basic idea of our approach  how to solve equation \eqref{1} generally on interval $[t_0, T]$, $T>t_0$ is the following:\\
If $ T \in (t_0, t_0+ \alpha]$ then we compute $u_1(t)$ as a desired approximate solution of \eqref{1} on $[t_0, T]$. \\
If not, we distinguish two cases, when the delays are commensurate and non-commensurate.

\subsection{Commensurate Delays}\label{commens}
If all delays $\tau_i$ are mutually commensurate, i.e. $\frac{\tau_i}{\tau_j}$ is a rational number for all $i,j=1, \ldots, r$, we define $\alpha_* = GCD (\tau_1, \ldots, \tau_r)$, where $GCD(\cdot)$ denotes greatest common divisor in the sense that $\alpha_*$ is the largest rational number for which there are positive integers $k_1, \ldots, k_r \in \mathbb{N}$ such that $\tau_i = k_i \cdot \alpha_*$ for $i = 1, \ldots, r$. Obviously $\alpha_* \leq \alpha$.

Further, we need to find $K \in \mathbb N$ such that $T \in (t_0+ K\alpha_*, t_0+ (K+1)\alpha_*]$.

First step is to find a solution $u_1$  in the interval $[t_0, t_0+ \alpha_*]$ in a way described above.

Next step is to find a solution $u_2$ in the interval $[ t_0+ \alpha_*, t_0 + 2\alpha_*]$. Then we solve the equation
\begin{equation}\label{7a}
u_2^{(n)}(t) = f(t,u_2(t),u_2'(t),\dots,u_2^{(n-1)}(t),\mathbf{u}_{1,1}(t-\tau_1),\mathbf{u}_{1,2}(t-\tau_2), \dots, \mathbf{u}_{1,r} (t-\tau_r) ),
\end{equation}
where either $\mathbf{u}_{1,i}(t-\tau_i )= (u_1(t-\tau_i),u_1'(t-\tau_i),\dots,u_1^{(m_i)}(t-\tau_i))$ or $\mathbf{u}_{1,i}(t-\tau_i )= (\phi(t-\tau_i),\phi'(t-\tau_i),\dots,\phi^{(m_i)}(t-\tau_i))$, $m_i \leq n$, $i=1,2,\dots,r$ depending on the situation: If the difference $t-\tau_i$ falls into $(t_0, t_0+ \alpha_*]$ for $t \in (t_0+ \alpha_*, t_0 + 2\alpha_*]$, the first expression is used. If $t-\tau_i \in [t_0-t^*,t_0]$ for $t \in (t_0+ \alpha_*, t_0 + 2\alpha_*]$, the second one is applied.

Applying DTM  to  equation \eqref{7a} leads to recurrence equation
\begin{equation}\label{eq8}
\frac{(k+n)!}{k!} U_2(k+n) = \mathcal F_2 \Bigl( k, U_2(k), U_2(k+1), \dots, U_2(k+n-1)\Bigr),
\end{equation}
where
$$
U_2(0)= u_1(t_0+\alpha_*)= u_{2,0},\ U_2(1)= u'_1(t_0+\alpha_*) =u_{2,1}, \dots , U_2(n-1)= \frac{1}{(n-1)!} u_1^{(n-1)}(t_0+\alpha_*) = \frac{u_{2,n-1}}{(n-1)!}.
$$

Using \eqref{eq8} we obtain a solution of \eqref{7a} in  the form
$$
u_2(t) = \sum_{k=0}^{\infty} U_2(k)\left(t-(t_0+\alpha_*)\right)^k, \quad t \in [t_0+ \alpha_*, t_0+ 2\alpha_*].
$$
Generally, the algorithm has the following form in $(j+1)$-th interval $[t_0+ j\alpha_*, t_0+ (j+1)\alpha_*]$ ($j =1, \ldots, K$):
We solve equation
\begin{equation}\label{eq7a}
u_{j+1}^{(n)}(t) = f(t,u_{j+1}(t),u'_{j+1}(t),\dots,u_{j+1}^{(n-1)}(t),\mathbf{u}_{j,1}(t-\tau_1),\mathbf{u}_{j,2}(t-\tau_2), \dots, \mathbf{u}_{j,r} (t-\tau_r) ),
\end{equation}
where 
\begin{equation}\label{eq7b}
\mathbf{u}_{j,i}(t-\tau_i )= (u_{l+1}(t-\tau_i),u'_{l+1}(t-\tau_i),\dots,u_{l+1}^{(m_i)}(t-\tau_i)) \quad
\text{if } t-\tau_i \in (t_0+ l\alpha_*, t_0+ (l+1)\alpha_*)
\end{equation}
for $t \in (t_0+ j\alpha_*, t_0+ (j+1)\alpha_*)$, $l \in \{ 0, \dots, j-1 \}$.

If $t-\tau_i \in [t_0-t^*,t_0]$ for $t \in (t_0+ j\alpha_*, t_0+ (j+1)\alpha_*]$, then again
$$
\mathbf{u}_{j,i}(t-\tau_i )=(\phi(t-\tau_i),\phi'(t-\tau_i),\dots,\phi^{(m_i)}(t-\tau_i)).
$$

Applying DTM  to  equation \eqref{eq7a} results in recurrence equation
\begin{equation}\label{eq9}
\frac{(k+n)!}{k!} U_{j+1}(k+n) = \mathcal F_{j+1} \Bigl( k, U_{j+1}(k), U_{j+1}(k+1), \dots, U_{j+1}(k+n-1)\Bigr),
\end{equation}
where
$$
U_{j+1}(0)= u_j(t_0+j\alpha_*)= u_{j+1,0},\ U_{j+1}(1)= u'_j(t_0+j\alpha_*) =u_{j+1,1}, \dots , U_{j+1}(n-1)= \frac{1}{(n-1)!} u_j^{(n-1)}(t_0+j\alpha_*) = \frac{u_{j+1,n-1}}{(n-1)!}.
$$

Using the above described algorithm, we get a solution of \eqref{eq7a}
$$
u_{j+1}(t) = \sum_{k=0}^{\infty} U_{j+1}(k)\left(t-(t_0+j\alpha_*)\right)^k, \quad t \in [t_0+ j\alpha_*, t_0+ (j+1)\alpha_*]
$$
for all $j \in \{ 0, \ldots, K \}$, since $T \in (t_0+ K\alpha_*, t_0+ (K+1)\alpha_*]$. Then  we get an approximate solution of \eqref{1} on the interval $[t_0, T]$ in the form
$$
u(t) = \left\{\begin{array}{ll}
        u_1(t),\ & t \in [t_0, t_0+ \alpha_*], \\
        u_2(t),\ & t \in [t_0+ \alpha_*, t_0+ 2\alpha_*]\\
        \vdots \\
        u_{K+1}(t), \ & t \in [t_0+ K\alpha_*, T].
        \end{array} \right.
$$

{\bf Remark.} Regarding \eqref{eq7b}, it is good to point out that  $\mathbf{u}_{j,i}(t-\tau_i )$, $j \in \{ 1, \ldots, K \}$, $i \in \{ 1, \ldots, r \}$ will always be transformed in such a way that we will utilize coefficients which already were computed in previous steps. This is especially important in computer implementation of the algorithm. Verification of this fact is not difficult:
\begin{lemma}\label{lemma2}
Let $t-\tau_i \in (t_0+ l\alpha_*, t_0+ (l+1)\alpha_*]$ for $t \in (t_0+ j\alpha_*, t_0+ (j+1)\alpha_*]$, where $l \in \{ 0, \dots, j-1 \}$. Let $p \in \mathbb N$. Then 
$$
\mathcal D \{u_{l+1}^{(p)} (t-\tau_i) \} [t_0+j\alpha_*] =\left\{ \frac{(k+p)!}{k!}U_{l+1}(k+p) \right\}_{k=0}^{\infty}.
$$
\end{lemma}
{\bf Proof.}
Since $\alpha_* = GCD (\tau_1, \ldots, \tau_r)$, then necessarily $t_0+ l\alpha_* + \tau_i = t_0+ j\alpha_*$, hence $\tau_i = (j-l) \alpha_*$. Then
$$
u_{l+1}^{(p)} (t - \tau_i) = u_{l+1}^{(p)} (t - (j-l)\alpha_*) = \sum_{k=0}^{\infty} \frac{(k+p)!}{k!} U_{l+1}(k+p)\bigl( (t-(j-l)\alpha_*)-(t_0+l\alpha_*) \bigr)^k = \sum_{k=0}^{\infty} \frac{(k+p)!}{k!} U_{l+1}(k+p)(t-(t_0+j\alpha_*))^k.
$$
\qed

\begin{example}
Suppose that we want to solve the following (neutral) equation
\begin{equation}\label{e1}
u'(t) = f \Bigl( t, u(t), u(t-1), u'(t-\frac{3}{2}) \Bigr)
\end{equation}
on $[0,3]$ with given initial function $g$ such that $u(t) = g(t)$ on $[-\frac{3}{2},0]$. We have $t_0 =0$, $T=3$, and two delays, $\tau_1 = 1$ and $\tau_2 = \frac{3}{2}$. Its ratio is a rational number, hence the delays are commensurate and we can follow the algorithm. We calculate that $\alpha_* = \frac{1}{2}$ (with $k_1 = 2$ and $k_2 = 3$) and $T \in (\frac{5}{2}, 3] = (t_0 + 5 \cdot \frac{1}{2}, t_0 + 6 \cdot \frac{1}{2}]$, so we have $K=5$. We need to solve the following sequence of equations:
\begin{align*}
u'_1(t) &= f \Bigl( t, u_1(t), g(t-1), g'(t-\frac{3}{2}) \Bigr) \ \text{ on } [0,\frac{1}{2}], \ u_1 (0) = g(0), \\
u'_2(t) &= f \Bigl( t, u_2(t), g(t-1), g'(t-\frac{3}{2}) \Bigr) \ \text{ on } [ \frac{1}{2},1],\  u_2 (\frac{1}{2}) = u_1 (\frac{1}{2}), \\
u'_3(t) &= f \Bigl( t, u_3(t), u_1(t-1), g'(t-\frac{3}{2}) \Bigr) \ \text{ on } [ 1, \frac{3}{2}], \ u_3 (1) = u_2 (1), \\
u'_4(t) &= f \Bigl( t, u_4(t), u_2(t-1), u'_1(t-\frac{3}{2}) \Bigr) \ \text{ on } [ \frac{3}{2}, 2], \ u_4 (\frac{3}{2}) = u_3 (\frac{3}{2}),\\
u'_5(t) &= f \Bigl( t, u_5(t), u_3(t-1), u'_2(t-\frac{3}{2}) \Bigr) \ \text{ on } [ 2, \frac{5}{2}], \ u_5 (2) = u_4 (2), \\
u'_6(t) &= f \Bigl( t, u_6(t), u_4(t-1), u'_3(t-\frac{3}{2}) \Bigr) \ \text{ on } [ \frac{5}{2}, 3], \ u_6 (\frac{5}{2}) = u_5 (\frac{5}{2}).
\end{align*}
Applying differential transformation to this sequence, we obtain a sequence of recurrence relations that we can solve using computer:
\begin{align*}
(k+1)U_1(k+1) &= \mathcal{F}_1 \Bigl( k, U_1(k), G_1(k), H_1(k) \Bigr) \ \text{ on } [0,\frac{1}{2}], \ U_1 (0) = g(0), \\
(k+1)U_2(k+1) &= \mathcal{F}_2 \Bigl( k, U_2(k), G_2(k), H_2(k) \Bigr) \ \text{ on } [ \frac{1}{2},1],\  U_2 (0) = u_1 (\frac{1}{2}), \\
(k+1)U_3(k+1) &= \mathcal{F}_3 \Bigl( k, U_3(k), U_1(k), H_3(k) \Bigr) \ \text{ on } [ 1, \frac{3}{2}], \ U_3 (0) = u_2 (1), \\
(k+1)U_4(k+1) &= \mathcal{F}_4 \Bigl( k, U_4(k), U_2(k), (k+1) U_1(k+1) \Bigr) \ \text{ on } [ \frac{3}{2}, 2], \ U_4 (0) = u_3 (\frac{3}{2}),\\
(k+1)U_5(k+1) &= \mathcal{F}_5 \Bigl( k, U_5(k), U_3(k), (k+1) U_2(k+1) \Bigr) \ \text{ on } [ 2, \frac{5}{2}], \ U_5 (0) = u_4 (2), \\
(k+1)U_6(k+1) &= \mathcal{F}_6 \Bigl( k, U_6(k), U_4(k), (k+1) U_3(k+1) \Bigr) \ \text{ on } [ \frac{5}{2}, 3], \ U_6 (0) = u_5 (\frac{5}{2}),
\end{align*}
where $G_1$ and $G_2$ are transformations of $k$-th derivative of $g(t-1)$ at $0$ and $\frac{1}{2}$ respectively, and $H_1$, $H_2$ and $H_3$ are transformations of $k$-th derivative of $g'(t-\frac{3}{2})$ at $0$, $\frac{1}{2}$ and $1$.
\end{example}

\subsection{Non-commensurate Delays}\label{noncommens}
If there are at least two delays $\tau_i$ and $\tau_j$ that are not commensurate, we cannot follow the algorithm described in \ref{commens}. We must proceed as follows:

First of all we need to find all (ordered) $r$-tuples $(K_1, \ldots, K_r) \in \mathbb{N}_0^r$ such that $t_0 +K_1 \tau_1 + \ldots + K_r \tau_r < T$. For all such $r$-tuples denote $\sigma_{(K_1, \ldots, K_r)} = t_0 + \sum\limits_{i=1}^{r} K_i \tau_i$, arrange these $\sigma_{(K_1, \ldots, K_r)}$ in a sequence having ascending order with respect to standard relation $\leq$ (if multiple $\sigma_{(K_1, \ldots, K_r)}$ have the same value, we keep only one) and relabel them, so that we obtain an ascending sequence $t_0 = \sigma_0 <\sigma_1<\sigma_2< \ldots <\sigma_j < \ldots <\sigma_K < T$, where $K \in \mathbb{N}$ depends on delays $\tau_1, \ldots, \tau_r$.

Denote $I_j = (\sigma_{j-1}, \sigma_j ]$ for $j=1, \ldots, K$, and put $I_{K+1} = (\sigma_K, T]$. Then we are looking for solution $u(t)$ in the form
$$
u(t) = \begin{cases}
       \ u_1(t),\ & t \in I_1, \\
       \ u_2(t),\ & t \in I_2, \\
       \quad \vdots \\
       \ u_{K+1}(t), \ & t \in I_{K+1},
        \end{cases}
$$
where solution $u_j$ in the $j$-th interval $I_j$ is obtained in the following way:

We solve equation
\begin{equation}\label{eq10}
u_{j}^{(n)}(t) = f(t,u_{j}(t),u'_{j}(t),\dots,u_{j}^{(n-1)}(t),\mathbf{u}_{j,1}(t-\tau_1),\mathbf{u}_{j,2}(t-\tau_2), \dots, \mathbf{u}_{j,r} (t-\tau_r) ),
\end{equation}
where 
\begin{equation}\label{eq11}
\mathbf{u}_{j,i}(t-\tau_i )= (u_{l}(t-\tau_i),u'_{l}(t-\tau_i),\dots,u_{l}^{(m_i)}(t-\tau_i)) \quad
\text{if } t-\tau_i \in I_l
\end{equation}
for $t \in I_j$, $l \in \{ 1, \dots, j-1 \}$, $j \in \{ 1, \ldots, K+1 \}$.

If $t-\tau_i \in [t_0-t^*,t_0]$ for $t \in I_j$, then again
$$
\mathbf{u}_{j,i}(t-\tau_i )=(\phi(t-\tau_i),\phi'(t-\tau_i),\dots,\phi^{(m_i)}(t-\tau_i)).
$$

Applying DTM  to  equation \eqref{eq10} results in recurrence equation
\begin{equation}\label{eq12}
\frac{(k+n)!}{k!} U_{j}(k+n) = \mathcal F_{j} \Bigl( k, U_{j}(k), U_{j}(k+1), \dots, U_{j}(k+n-1)\Bigr),
\end{equation}
where
$$
U_{j}(0)= u_{j-1}(\sigma_{j-1})= u_{j,0},\ U_{j}(1)= u'_{j-1}(\sigma_{j-1}) =u_{j,1}, \dots , U_{j}(n-1)= \frac{1}{(n-1)!} u_{j-1}^{(n-1)}(\sigma_{j-1}) = \frac{u_{j,n-1}}{(n-1)!},
$$
$j \in \{ 1, \ldots, K+1 \}$, and $u_0 (t) = \phi (t)$.

Using the above described algorithm, we get a solution of \eqref{eq10}
$$
u_{j}(t) = \sum_{k=0}^{\infty} U_{j}(k)\left(t-\sigma_{j-1})\right)^k, \quad t \in I_j
$$
for all $j \in \{ 1, \ldots, K+1 \}$.

{\bf Remark.} At this point, it is important to notice that with \eqref{eq11} we are in a different situation than in \eqref{eq7b}. We will still utilize coefficients computed in previous steps for transforming $\mathbf{u}_{j,i}(t-\tau_i )$, $j \in \{ 1, \ldots, K \}$, $i \in \{ 1, \ldots, r \}$. However, the relation is different from formula derived for \eqref{eq7b} in Lemma \ref{lemma2}. We formulate it in the following lemma:
\begin{lemma}\label{lemma3}
Let $t-\tau_i \in I_l$ for $t \in I_j$, where $l \in \{ 1, \dots, j-1 \}$. Let $p \in \mathbb N$. Then 
\begin{equation}\label{l3_1}
\mathcal D \{u_{l}^{(p)} (t-\tau_i) \} [\sigma_{j-1}] =\left\{ \frac{(k+p)!}{k!}\sum\limits_{y=k+p}^{\infty} \binom{y}{k+p} (\sigma_{j-1} - \sigma_{l-1} - \tau_i)^{y-k-p} U_{l}(y) \right\}_{k=0}^{\infty}.
\end{equation}
\end{lemma}
{\bf Proof.}
We have
$$
u_{l} (t-\tau_i) = \sum\limits_{k=0}^{\infty} U_l (k) (t-\tau_i-\sigma_{l-1})^{k} = \sum\limits_{k=0}^{\infty} U_l (k) (t-\sigma_{j-1} +\sigma_{j-1} - \sigma_{l-1} -\tau_i)^k = \sum\limits_{k=0}^{\infty} U_l (k) \sum\limits_{y=0}^{k} \binom{k}{y} (t-\sigma_{j-1})^y (\sigma_{j-1} - \sigma_{l-1} -\tau_i)^{k-y}.
$$
Rearranging the order of summation and relabeling ($y \rightarrow k$, $k \rightarrow y$), we obtain
$$
u_{l} (t-\tau_i) = \sum\limits_{k=0}^{\infty} \left( \sum\limits_{y=k}^{\infty} \binom{y}{k} (\sigma_{j-1} - \sigma_{l-1} -\tau_i)^{y-k} U_l (y) \right) (t-\sigma_{j-1})^k.
$$
Now \eqref{l3_1} is a consequence of formula $i)$ from Lemma \ref{lemma1}.
\qed

{\bf Remark.} It may happen that $\sigma_{j-1} - \sigma_{l-1} - \tau_i$ vanishes in some cases. Under such circumstances, the formula \eqref{l3_1} is simplified to
$$
\mathcal D \{u_{l}^{(p)} (t-\tau_i) \} [\sigma_{j-1}] =\left\{ \frac{(k+p)!}{k!} U_l (k+p) \right\}_{k=0}^{\infty}.
$$

{\bf Remark.} In practical applications, upper bound of the sum in \eqref{l3_1} must be finite, not exceeding $N$. It means that, in contrast to \eqref{eq7b}, the result values of transformation procedure of  $\mathbf{u}_{j,i}(t-\tau_i )$ are only approximations. However, as we already work with irrational numbers (some delays are not commensurate), this is only a minor nonsignificant issue.

\begin{example}
Suppose that we want to solve the following equation
\begin{equation}\label{e2}
u'(t) = f \Bigl( t, u(t), u(t-3), u(t-\pi) \Bigr)
\end{equation}
on $[0,9]$ with given initial function $g$ such that $u(t) = g(t)$ on $[-\pi,0]$.  We have $t_0 =0$, $T=9$, and two delays, $\tau_1 = 3$ and $\tau_2 = \pi$. Its ratio is an irrational number, hence the delays are not commensurate and we cannot follow the algorithm described in Subsection \ref{commens}. Instead of it we must use the ideas introduced in Subsection \ref{noncommens}.

First we find the $2$-tuples (ordered pairs): $(0,0)$, $(1,0)$, $(0,1)$, $(2,0)$, $(1,1)$, $(0,2)$, which give us the following sequence: $\sigma_{(0,0)}=\sigma_0 = t_0 = 0$, $\sigma_{(1,0)}=\sigma_1 = 3$, $\sigma_{(0,1)}=\sigma_2 = \pi$, $\sigma_{(2,0)}=\sigma_3 = 6$, $\sigma_{(1,1)}=\sigma_4 = 3 + \pi$, $\sigma_{(0,2)}=\sigma_5 = 2 \pi$. We see that $K=5$.

Now, we have $I_1 = (0,3]$, $I_2 = (3,\pi]$, $I_3 = (\pi,6]$, $I_4 = (6,3+\pi]$, $I_5 = (3+\pi,2\pi]$, $I_6 = (2\pi,9]$. We need to solve the following sequence of equations:
\begin{align*}
u'_1(t) &= f \Bigl( t, u_1(t), g(t-3), g(t-\pi)) \Bigr) \ \text{ on } I_1, \ u_1 (0) = g(0), \\
u'_2(t) &= f \Bigl( t, u_2(t), u_1(t-3), g(t-\pi)) \Bigr) \ \text{ on } I_2,\  u_2 (3) = u_1 (3), \\
u'_3(t) &= f \Bigl( t, u_3(t), u_1(t-3), u_1(t-\pi)) \Bigr) \ \text{ on } I_3, \ u_3 (\pi) = u_2 (\pi), \\
u'_4(t) &= f \Bigl( t, u_4(t), u_2(t-3), u_1(t-\pi)) \Bigr) \ \text{ on } I_4, \ u_4 (6) = u_3 (6),\\
u'_5(t) &= f \Bigl( t, u_5(t), u_3(t-3), u_2(t-\pi)) \Bigr) \ \text{ on } I_5, \ u_5 (3+\pi) = u_4 (3+\pi), \\
u'_6(t) &= f \Bigl( t, u_6(t), u_3(t-3), u_3(t-\pi)) \Bigr) \ \text{ on } I_6, \ u_6 (2\pi) = u_5 (2\pi).
\end{align*}
\end{example}

\section{Test Examples}
In this section we will apply the derived algorithms to delayed and neutral logistic differential equations to show the efficiency of our method in comparison with another numerical solution methods such as variational  iteration method, polynomial least squares method, Adomian decomposition method and solutions of \eqref{l2} and \eqref{l5} using software Matlab, package DDE23.

\begin{example}\label{ex3}
First we turn back to Hutchinson's logistic delay population model
\begin{equation}\label{p1}
u'(t) = r u(t) \left( 1- \frac {u(t-\tau)}{K}\right), \quad t \in [0,T] \tag{\ref{l2}}
\end{equation}
subject to the initial function
$$
u(t) =\phi(t), \quad t \in [-\tau, 0].
$$
In general, differential transformation of equation \eqref{p1} in combination with method of steps give
\begin{equation}\label{r1}
U(k+1) = \frac{1}{k+1}\left[ r U(k) - \frac{r}{K} \sum_{l=0}^k F(l)U(k-l)\right], \ k \geq 0
\end{equation}
where $F(l)$ is differential transformation of the function $\phi(t-\tau)$. In the following we consider $K=0.5, r=2, \tau = 0.1, \phi(t)=1 $ as Dehghan and Salehi \cite{dehghan} or Caruntu and Bota \cite{caruntu}
who solved the given problem using Adomian decomposition metod, variational iteration method and the polynomial least squares method, respectively.

In this case, equation \eqref{p1} becomes

\begin{equation}\label{p2}
u'(t) = u(t) \left(2-4u(t-0.1)\right)
\end{equation}
subject to intial function
$$
u(t) = \phi(t) =1, \ t \in [-0.1,0].
$$
As we have only one constant delay, we can follow the algorithm described in Subsection \ref{commens} (commensurate delays).
Applying the method of steps and the differential transformation to equation \eqref{p2}, we have
\begin{equation}\label{d1}
U_1(k+1)= \frac{-2}{k+1}U_1(k), \ k \geq 0.
\end{equation}
We will see that it is enough to  calculate the truncated series solution of equation \eqref{p2} for $N=3$ to obtain simplier and more accurate results than using the above mentioned methods.\\
From \eqref{d1} we obtain 
 $$U_1(0) = 1,\ U_1(1) = -2,\ U_1(2)= 2,\ U_1(3)=-4/3.$$
Then
$$
u_1(t) = \sum_{k=0}^{3} U_1(k) \cdot t^k = 1-2t+2t^2 -\frac{4}{3}t^3 + O(t^4), \quad t \in [0,0.1].
$$
Following the algorithm, we will solve  equation \eqref{p2} on the interval $[0.1,0.2]$ with the prescribed constants and the initial function
$$
\phi(t) = u_1(t) =1-2t+2t^2 -\frac{4}{3}t^3, \quad t \in [0, 0.1].
$$
Then  equation \eqref{p2} has the form
$$
u_2'(t) =  u_2(t) \left[ -2+8(t-0.1) -8(t-0.1)^2 +16/3(t-0.1)^3\right]
$$ 
and 
$$
U_2(k+1) =   \frac{1}{k+1} \left[ \sum_{l=0}^k\left(-2\delta(l)+8\delta(l-1)-8\delta(l-2)+\frac{16}{3}\delta(l-3)\right)U_2(k-l)\right],\  k \geq 0.
$$
Then
\begin{align*}
U_2(0)&=0.81867,  \\
U_2(1)&= -2U_2(0)= -1.63734, \\  
U_2(2)&= -U_2(1)+4U_2(0)= 4.91202,  \\
U_2(3)&= -\frac{2}{3}U_2(2)+\frac{8}{3}U_2(1)-\frac{8}{3}U_2(0) = -9.82404,  \\
U_2(4)& = 1/4\left( -2U_2(3) +8U_2(2) -8U_2(1) +16/3U_2(0)\right)= 19.10230, \\
U_2(k)&= 0 \ {\rm for} \ k \ \geq 5 \\
\end{align*}
and
$$
u_2(t) = \sum_{k=0}^4 U_2(k) \cdot (t-0.1)^k = 0.81867 -1.63734(t-0.1)+4.91202(t-0.1)^2 -9.82404(t-0.1)^3 +19.10230(t-0.1)^4
$$
for $ t\in [0.1,0.2]$.\\
Now we will solve equation \eqref{p2} with the initial function
$$
\phi(t) = u_2(t) = 0.81867 -1.63734(t-0.1)+4.91202(t-0.1)^2 -9.82404(t-0.1)^3 +19.10230(t-0.1)^4, \ t \in [0.1,0.2].
$$
Then equation \eqref{p2} has the form
$$
u_3'(t) =  u_3(t) \left[ -1.27468+6.54936(t-0.2)-19.64808(t-0.2)^2+39.29616(t-0.2)^3-76.40920(t-4)^4 \right]
$$ 
and
\begin{align}\label{u1}
U_3(k+1) &=   \frac{1}{k+1}  \left[ \sum_{l=0}^k\left( -1.27468\delta(l)+6.54936\delta(l-1)-19.64808\delta(l-2)+39.29616\delta(l-3) \right.\right. \\
& \left.\left.   -76.40920\delta(l-4)\right)U_2(k-l)\right] \nonumber
\end{align}
for $k \geq 0 $. From \eqref{u1} and the initial condition we obtain
\begin{align*}
U_3(0)&=0.69614,\ U_3(1)=-0.88736,\ U_3(2)=2.84506,\ U_3(3)=-7.70533,\ U_3(4)=18.31142,\\
 U_3(5)&=-43.55350, \ U_3(k)=0\ {\rm for} \ k \geq 6.
\end{align*}
Then
\begin{align*}
u_3(t) &= \sum_{k=0}^5 U_3(k) \cdot (t-0.2)^k= 0.69614-.88736(t-0.2)+2.84506(t-0.2)^2-7.70533(t-0.2)^3\\
           &+18.31142(t-0.2)^4-43.55350(t-0.2)^5,\ t \in [0.2, 0.3].
\end{align*}
Repeting the previous procedure with initial function
$$
\phi(t)= u_3(t), \ t \in [0.2,0.3]
$$
we get
\begin{align*}
u_4(t) & = \sum_{k=0}^6 U_4(k) \cdot (t-0.3)^k=0.62954-0.49580(t-0.3)+1.31175(t-0.3)^2-3.31776(t-0.3)^3+8.076134(t-0.3)^4\\
&-18.88654(t-0.3)^5 +44.61284(t-0.3)^6
\end{align*}
for $ t\in [0.3,0.4].$  Generally for intial function $ \phi(t) = u_{k-1}(t)$ for $ t\in [(k-2)0.1,(k-1)0.1] $ we obtain truncated series solution 
$$
u_k(t) = \sum_{j=0}^{k+2} U_k(j) \cdot (t-(k-1)0.1)^j, \ t \in [(k-1) \cdot 0.1, k \cdot 0.1], \ k \geq 2.
$$
Comparison of results obtained by differential tranformation method (DTM), by other methods (polynomial least square method - PLSM, Adomian decomposition method - ADM, variational iteration method - VIM) and by Matlab module DDE23 is done in Table 1.
\begin{center}
Table 1.\\
\begin{tabular}{|c|c|c|c|c|c|}
\hline t    &   DTM        & PLSM        & ADM     & VIM      & Matlab DDE23  \\ 
\hline\hline 0.0 &1.00000  &1.00000&     1.00000&  1.00000&  1.00000  \\
\hline 0.05  &0.90484&     0.84150&     0.79988&   0.63797&  0.90484 \\
\hline 0.10  &0.81867&     0.74007&     0.67999&   0.58790&  0.81873 \\
\hline 0.15  &0.74797&     0.67671&     0.61608&   0.55611&  0.74802\\
\hline 0.20  &0.69614&     0.63760&     0.56954&   0.53635&  0.69591 \\
\hline 0.25  &0.65802&     0.61308&     0.53249&   0.52433&  0.65782\\
\hline 0.30  &0.62954&     0.59668&     0.50241&   0.51711&  0.62940\\
\hline 0.35 &0.60765 &     0.58437&     0.47988&   0.51276&  0.60765 \\
\hline 0.40 &0.59042 &     0.57385&     0.46124&   0.51001&  0.59052\\
\hline 0.45 &0.57663 &     0.56405&     0.45024&   0.50808&  0.57673\\
\hline 0.50 &0.56536 &     0.55467&     0.44803&   0.50651&  0.56546
\\ \hline
\end{tabular}
\end{center}
The solutions are plotted in Figure 1.

Looking at the table and the figure, we can observe that the proposed method gives results close to the solution produced by Matlab DDE23 whereas ADM and VIM are not accurate enough and especially ADM is not suitable for solving equation \eqref{p2} at all. Moreover, already the first and the second iteration formulas of a solution of \eqref{p2} using ADM and VIM are very complicated (see \cite{dehghan}).
\begin{figure}
\label{Fig1}
\centering
\includegraphics[bb=300 230 300 530]{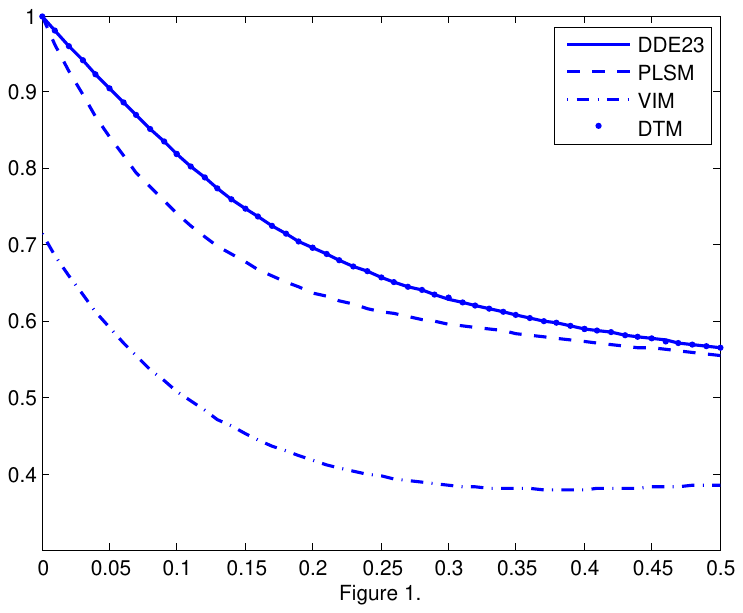}
\end{figure}
\end{example}

\begin{example}
As an example of neutral type logistic equation, we have chosen equation studied by Gy\"ori and Ladas \cite{gyori2}
\begin{equation}\label{p2.1}
u'(t) = u(t)\left\{ r \left[1-\frac{u(t-\sigma)}{K} \right]+ c \frac{u'(t-\tau)}{u(t-\tau)}\right\}, \tag{\ref{l5}}
\end{equation}
with  initial function
$$
u(t) =\phi(t), \quad t \in [-\gamma, 0], \quad \gamma = \max \{ \sigma, \tau \}.
$$
Oscillatory and bifurcation properties of solutions of this initial problem were recently studied by Wang and Wei \cite{wang1} and Wang \cite{wang2}. In general case of initial function $\phi (t)$ depending on $t$, differential transformation of equation \eqref{p2.1} would be more complicated than differential transformation of equation \eqref{p1} in Example \ref{ex3}, since there is $u(t-\tau)$ in the denominator. Therefore we will proceed with the same conditions as Wang \cite{wang2}: $K=3$, $r=0.45$, $\sigma = 2$, $c=0.3$, $\tau = 1$ and $\phi(t)=2.3$ for $t \in [-2,0]$. Let $T=2$.

 In this case, equation \eqref{p2.1} has the form
\begin{equation}\label{p2.2}
u'(t) = u(t)\left\{ 0.45 \left[1-\frac{u(t-2)}{3} \right]+ 0.3 \frac{u'(t-1)}{u(t-1)}\right\},
\end{equation}
subject to initial function
$$
u(t) =\phi(t) = 2.3, \quad t \in [-2, 0].
$$
As we have two commensurate constant delays $\sigma =\tau_1 =2$ and $\tau =\tau_2 =1$, we can follow the algorithm described in Subsection \ref{commens}. We have $t_0 =0$, $T=2$, we calculate that $\alpha_* =1$ (with $k_1 =2$ and $k_2 =1$) and $T \in (1,2] = (t_0 + 1 \cdot 1, t_0 + 2 \cdot 1]$, so we have $K=1$.

Applying the method of steps and the differential transformation to equation \eqref{p2.2}, we have
\begin{equation}\label{p2.3}
U_1(k+1)= \frac{0.105}{k+1}U_1(k), \ k \geq 0.
\end{equation}
From \eqref{p2.3} and intial condition $U_1 (0)=2.3$ we obtain
$$
U_1(0) = 2.3,\ U_1(1) = 2.3 \cdot 0.105,\ U_1(2)= 2.3 \cdot \frac{0.105^2}{2},\ U_1(3)=2.3 \cdot \frac{0.105^3}{3!}, \ U_1(4)=2.3 \cdot \frac{0.105^4}{4!},\ldots
$$
Here we can even see the closed form (analytical) solution:
$$
u_1 (t) = \sum_{k=0}^{\infty} U_1 (k) \cdot t^k = \sum_{k=0}^{\infty} 2.3 \cdot \frac{(0.105 \cdot t)^k}{k!} =2.3 \cdot \e^{0.105 \cdot t}, \quad t \in (0,1].
$$
Since $K=1$, we need to calculate $u_2$ too. Following the algorithm, we will solve equation \eqref{p2.2} on interval $(1,2]$ with the prescribed constants and initial functions
$$
u(t) =
\begin{cases}
\phi(t) = 2.3, \quad t \in [-2, 0] \\
u_1 (t) = 2.3 \cdot \e^{0.105 \cdot t}, \quad t \in (0,1].
\end{cases}
$$
Equation \eqref{p2.2} then has the form
$$
u_2'(t) = u_2 (t)\left\{ 0.45 \left[1-\frac{2.3}{3} \right]+ 0.3 \frac{0.105 \e^{0.105 \cdot (t-1)}}{\e^{0.105 \cdot (t-1)}} \right\} = 1.3 \cdot 0.105 \cdot u_2 (t)
$$
and
$$
U_2 (k+1) = \frac{1.3 \cdot 0.105}{k+1} U_2 (k), \ k \geq 0.
$$
We calculate
\begin{align*}
U_2(0) &= u_1 (1) = 2.3 \cdot \e^{0.105},\ U_2(1) = 2.3 \cdot \e^{0.105} \cdot 1.3 \cdot 0.105,\ U_2(2)= 2.3 \cdot \e^{0.105} \cdot \frac{(1.3 \cdot 0.105)^2}{2},\\
U_2(3)&=2.3 \cdot \e^{0.105} \cdot \frac{(1.3 \cdot 0.105)^3}{3!},
U_2(4) =2.3 \cdot \e^{0.105} \cdot \frac{(1.3 \cdot 0.105)^4}{4!},\ldots
\end{align*}
Again, we can guess the closed form solution:
$$
u_2 (t) = \sum_{k=0}^{\infty} U_2 (k) \cdot (t-1)^k = \sum_{k=0}^{\infty} 2.3 \cdot \e^{0.105} \cdot \frac{\bigl( 1.3 \cdot 0.105 \cdot (t-1) \bigr)^k}{k!} =2.3 \cdot \e^{0.105} \cdot \e^{1.3 \cdot 0.105 \cdot (t-1)}, \quad t \in (1,2].
$$
For computer simulations, we transcribed the algorithm in language of computer algebra system Maple using only cycles and \emph{PolynomialTools} package. We made a comparison with analytical solution obtained using procedure "dsolve" in Maple and with numerical solution  produced by Matlab DDE23. The comparison is shown in Table 2.
\begin{center}
Table 2.\\
\begin{tabular}{|c|c|c|c|c|c|}
\hline t      & DTM   &      Maple &     Matlab DDE23  \\ 
\hline
\hline 0.0  &2.3000&     2.3000&     2.3000 \\
\hline 0.2  &2.3488&     2.3488&     2.3488 \\
\hline 0.4  &2.3987&     2.3987&     2.3987 \\
\hline 0.6  &2.4496&     2.4496&     2.4496 \\
\hline 0.8  &2.5015&     2.5015&     2.5015 \\
\hline 1.0  &2.5546&     2.5546&     2.5546 \\
\hline 1.2  &2.6253&     2.6253&     2.6251 \\
\hline 1.4  &2.6980&     2.6980&     2.6977 \\
\hline 1.6  &2.7727&     2.7727&     2.7724 \\
\hline 1.8  &2.8494&     2.8494&     2.8492 \\
\hline 2.0  &2.9283&     2.9283&     2.9278
\\ \hline
\end{tabular}
\end{center}
Looking at the table, we can observe that the algorithm gives results close to the numerical solution produced by Matlab DDE23 and identical with the results obtained by analytical tools in Maple. 
\end{example}

{\bf Remark.} Before we proceed with a conclusion, it is good to mention recent works concerning discrete equivalent of continuous logistic equation without delay (\cite{petro1}, \cite{petro2}). Although the authors used more general approach involving Hilbert spaces, the derived transformation rules are identical to DTM formulas. Since DTM combined with method of steps is applicable for functional differential equations, our suggestion is that this general approach involving Hilbert spaces could be extended to delayed or neutral logistic equation too.

\section{Conclusion}
\begin{itemize}
\item
The approach presented in this paper is powerful and efficient semi-analytical technique generally convenient for numerical approximation of a unique solution of initial problem for a class of functional differential equations including delayed and neutral differential equations.
\item
Two algorithms were developed for models with several constant delays. We distinguished two cases when the delays are commensurate or non-commensurate.
\item
Using presented algorithms, we are able to obtain approximate solution, but sometimes there even is a possibility to identify unique solution of initial problem in closed form.
\item
We do not need to calculate multiple integrals or derivatives, hence less computational work is demanded compared to other popular methods (variational iteration method, homotopy perturbation method, homotopy analysis method, Adomian decomposition method).
\item
Compared to any purely numerical method, a specific advantage of this technique is that the solution is always a piecewise smooth function.
\item
Further investigation can be focused on development of the presented technique for equation \eqref{1} with distributed, state dependent or time dependent delays.
\end{itemize}

\section{Acknowledgements}
\label{sec6}
This research was carried out under the project CEITEC 2020 (LQ1601) with financial support from the Ministry of Education, Youth and Sports of the Czech Republic under the National Sustainability Programme II.
This support is gratefully acknowledged.





\bibliographystyle{elsarticle-num}
\bibliography{<your-bib-database>}



\end{document}